\documentclass[12pt]{amsart}
\usepackage{amsfonts}
\usepackage{amssymb}

\usepackage{amsmath}
\usepackage{epsf}
\usepackage{epsfig}
\usepackage{colordvi}
\usepackage{fancybox}

\newcommand{\R}{\mathbb R}
\newcommand{\C}{\mathbb C}
\newcommand{\D}{\mathbb D}

\newcommand{\jj}{\mathrm{j}}

\newtheorem{thm}{Theorem}[section]
\newtheorem{cor}[thm]{Corollary}

\newtheorem{prop}[thm]{Proposition}
\theoremstyle{definition}

\theoremstyle{remark}
\newtheorem{remark}{Remark}[section]

\newcommand{\ds}{\displaystyle}

\begin{document}

\title[ Minimal Timelike Surfaces in $\mathbb R^3_1$ and their Canonical Parameters]
{Minimal Timelike Surfaces in the Lorentz-Minkowski 3-space and their Canonical Parameters}%

\thanks{2020 {\it Mathematics Subject Classification}: 53A10; 53B30; 53C50}

\author{Ognian Kassabov and Velichka Milousheva}%

\address{Institute of Mathematics and Informatics, Bulgarian Academy of Sciences,
Acad. G. Bonchev Str. bl. 8, 1113, Sofia, Bulgaria}
\email{vmil@math.bas.bg}
\email{okassabov@math.bas.bg}

\keywords{Timelike surfaces, canonical  parameters, Weierstrass formula}%

\begin{abstract}
We study minimal timelike surfaces in  $\mathbb R^3_1$ using  a special Weierstrass-type formula in terms of holomorphic functions defined in the algebra of the double (split-complex) numbers.
We present a method of obtaining an equation of a minimal timelike surface in terms of canonical parameters, which play a role similar to the role of the natural parameters of curves  in $\mathbb R^3$.  
Having one holomorphic function that generates a minimal timelike surface,
we find all holomorphic functions that generate the same surface. In this way we give  a correspondence 
between a minimal timelike surface and a class of holomorphic functions. 
As an application, we prove that the Enneper surfaces are the only minimal 
timelike surfaces in  $\mathbb R^3_1$ with polynomial parametrization of degree 3 in isothermal parameters.

\end{abstract}
\maketitle

\section{Introduction}
 
The study of minimal surfaces is one of the main topics in classical differential geometry which  goes back to the 18th century. Lagrange  initiated in 1760 the study of minimal surfaces in Euclidean 3-space and found the minimal surface equation  when he looked  for a necessary condition for minimizing the area functional. 
He showed that a minimal surface parametrized as a graphic $x=(u,v,\varphi(u,v))$ satisfies the following equation, known nowadays as the Lagrange's equation, 1762:
 $$(1+\varphi_v^2)\varphi_{uu}-2\varphi_u \varphi_v \varphi_{uv}+(1+\varphi_u^2)\varphi_{vv}=0.$$
The link between curvature and minimal surfaces was made by Meusnier in 1776 who proved that the  Lagrange's equation implies that the mean curvature is zero everywhere on a minimal surface. 
Usually, minimal surfaces are  defined as surfaces with zero mean curvature, but they are also characterized as surfaces of minimal surface area for given boundary conditions, as a critical point of the area functional, or as a graphic of the solution of a differential equation. 

The Weierstrass representation formula (1866) describes minimal surfaces in terms of two holomorphic functions $f(z)$ and $g(z)$ as follows \cite{Weier}:
$$
	\Psi(z)=\Re \int_{z_0}^z \, 
	\left( \frac12 f(z)(1-g^2(z)),\frac i2f(z)(1+g^2(z)),f(z)g(z) \right) \, dz.
$$

\vskip 1mm

The theory of minimal surfaces in real space forms have been attracting the attention of many mathematicians for more than two  centuries and have inspired many authors to study minimal surfaces in other ambient spaces. 
In  the last years, great attention is  paid to Lorentz surfaces in pseudo-Euclidean spaces, since pseudo-Riemannian geometry has many important applications in Physics, especially in problems related to General Relativity. 

However, the local geometry of surfaces in the Lorentz-Minkowski space $\mathbb R^3_1$  is much more complicated than that in the Euclidean space $\R^3$, since in $\R^3_1$ the vectors have different casual characters (spacelike, timelike or lightlike), which yield more  cases to be considered. One could consider spacelike, timelike or lightlike surfaces in $\R^3_1$. In the present paper we consider minimal timelike surfaces. Although in the timelike case the minimal surfaces neither maximize nor minimize surfaces area, they have many geometric properties similar to minimal surfaces in the Euclidean space $\R^3$. For example, Weierstrass representation formula was introduced by M. Magid in \cite{Mag}. While the classical Weierstrass representation theorem uses a relationship between holomorphic functions and solutions of certain elliptic PDEs, in the case of minimal timelike surfaces in $\R^3_1$ the problem of finding a  Weierstrass type representation is related to solving a certain system of hyperbolic PDEs. This shows that the difference between the pseudo-Riemannian and Riemannian case is as between hyperbolic and elliptic PDEs. In \cite{Mag}, M. Magid obtained local results. A global version of the  Weierstrass representation theorem for Lorentz surfaces in $\R^3_1$ is given by J. Konderak in \cite{Kond}. After the work of Konderak, many researches have studied spacelike and timelike surfaces in a 3\nobreakdash-dimensional Lorentz space via Weierstrass-type formulas, see e.g. \cite{Cintra-Onnis}, \cite{Lira-Melo-Merc}, \cite{Ship-Pack}. 

In the Euclidean space, the classical Bj\"orling problem, proposed by Bj\"orling in 1844,  is related to the construction of  a minimal surface in $\R^3$ containing a prescribed analytic strip.
The Bj\"orling problem  for timelike surfaces in the Lorentz-Minkowski spaces $\R^3_1$ was solved in \cite{Ch-D-Mag}, where a representation formula was obtained by use of split-complex numbers and natural split-complex extensions. The split-complex numbers are also known as Lorentz numbers, 
para-complex, double or hyperbolic numbers, and play a role similar to that played by the ordinary complex numbers in the spacelike case. The algebra of Lorentz numbers is often used  in the study of timelike surfaces, see for instance \cite{Erd}, \cite{Fu-In}, \cite{In-Tod}.

\vskip 2mm
In differential geometry, in the study of surfaces in both Euclidean and pseudo-Euclidean geometry, it is appropriate to use different types of special parameters: isothermal, principal, asymptotic, geodesic, etc., in different particular cases. In this paper, we deal with a special kind
of parameters, which play a role similar to the role of the natural parameters of a curve which are essential in the theory of curves in $\mathbb R^3$. Note that two important invariants --
the curvature and the torsion of a regular curve, parametrized by natural parameters,  determine the curve up a motion in the space. 

Actually,  natural parameters of a surface are not known in the general case, even in the Euclidean space $\R^3$. However, 
such parameters were introduced for the Weingarten surfaces which form a wide class of surfaces in  $\mathbb R^3$ \cite{G-M}. 
Some results in   \cite{G-M} were specialized in \cite{G-2008} for the class of minimal surfaces (i.e. surfaces with vanishing mean curvature)
in  $\mathbb R^3$ and in this case the special 
parameters are called canonical principal parameters. They are determined up to renumbering, sign and additive constants.  
The normal curvature, expressed in terms of canonical principal parameters, satisfies a special partial differential equation determining 
the surface up to a motion in  $\mathbb R^3$. Note that a surface may be parametrized in quite different ways and it is not
easy to say whether or not two parametric equations define one and the same surface. The normal curvature of minimal surfaces gives 
a simple method to answer this question: one may simply calculate the normal curvatures of two  minimal surfaces in terms of canonical 
principal parameters and then compare the results. The problem is that usually a minimal surface is defined in arbitrary
(not canonical) parameters. An effective method to solve this problem was proposed in \cite{OK}.
  
Similar questions may be considered for surfaces with vanishing mean curvature in the Lorentz-Minkowski space  $\mathbb R^3_1$.  An important tool to investigate the minimal 
surfaces in  $\mathbb R^3$ are the complex numbers and the complex functions (on the field of complex numbers),
see e.g. \cite{G-A-S}. These functions are also convenient
in the study of spacelike surfaces in  $\mathbb R^3_1$ with vanishing mean curvature -- the maximal surfaces.

In the present paper, we are interested in minimal timelike surfaces in  $\mathbb R^3_1$, so we apply the theory of functions on the algebra of double numbers. 
We use also the special Weierstrass formula of G. Ganchev, proposed in \cite{G-time-like}, where a new
approach to timelike surfaces is given and a determination of a minimal timelike surface via canonical
parameters is established. We find a method of obtaining an equation of a minimal timelike surface in canonical parameters. 

Having one holomorphic function (defined on a domain of
the plane of double numbers) that generates a minimal timelike surface,
we find all holomorphic functions that generate the same surface. Thus we obtain a correspondence 
between a minimal timelike surface and a class of holomorphic functions. 
As an application, we prove that the Enneper surfaces are the only minimal 
timelike surfaces in  $\mathbb R^3_1$ with polynomial parametrization of degree 3 in isothermal parameters.

\setcounter{equation}{0}
\section{Preliminaries}

We deal with the 3-dimensional Lorentz-Minkowski space  $\mathbb R^3_1$, endowed with the standard flat metric of signature  $(2,1)$:
$$
	\langle x, y \rangle = -x_1y_1+x_2y_2+x_3y_3 \ .
$$
The considerations in this paper are local and all functions are supposed to be of class $C^{\infty}$.

A regular surface $S$ in  $\mathbb R^3_1$ is said to be:

\hskip 6mm 
- \textit{timelike}, if the restriction of $\langle ., . \rangle$ to each tangent space of $S$ is indefinite;
 
\hskip 6mm 
- \textit{spacelike}, if the restriction of $\langle ., . \rangle$  to each tangent space of $S$ is positive definite;

\hskip 6mm 
- \textit{lightlike}, if the restriction of $\langle ., . \rangle$  to each tangent space of $S$ is degenerate.

In what follows, we suppose that the surface $S$ is timelike and is defined by the parametric equation
$$
	{\bf x=x}(u,v)=(x_1(u,v),x_2(u,v),x_3(u,v)) , \qquad 
	(u,v) \in U \subset \mathbb{R}^2 .
$$
We denote the derivatives of the vector function ${\bf x=x}(u,v)$ by
$$
	{\bf x}_u=\frac{\partial \bf x}{\partial u} ,  \qquad
	{\bf x}_v=\frac{\partial \bf x}{\partial v} ,  \qquad
	{\bf x}_{uv}=\frac{\partial^2 \bf x}{\partial u\partial v} ,  \ ...
$$
The coefficients of the first fundamental form are given by
$$
	E=\langle {\bf x}_u, {\bf x}_u \rangle, \qquad F= \langle {\bf x}_u, {\bf x}_v \rangle, \qquad G= \langle {\bf x}_v, {\bf x}_v \rangle .
$$
Denote by ${\bf U}$ the unit normal to the surface, i.e. 
$$
	{\bf U} = \frac{{\bf x}_u\times {\bf x}_v}{|{\bf x}_u\times {\bf x}_v|} .
$$
Then, the coefficients of the second fundamental form are
$$
	L= \langle {\bf U}, {\bf x}_{uu} \rangle, \qquad
	M= \langle {\bf U}, {\bf x}_{uv} \rangle, \qquad
	N= \langle {\bf U}, {\bf x}_{vv} \rangle .
$$
The Gauss curvarture and the mean curvature of $S$ are defined respectively by
$$
	K=\frac{LN-M^2}{EG-F^2}, \qquad
	H=\frac{EN-2FM+GL}{2(EG-F^2)} .
$$
The surface $S$ is said to be {\it minimal} if the mean curvature vanishes identically.

Probably, the most used parameters of a surface are the isothermal parameters. The coefficients of the first fundamental form of a 
timelike surface in $\mathbb R^3_1$ parametrized in terms of isothermal parameters satisfy $E=-G$, $F=0$. 
In the study of minimal timelike surfaces in $\mathbb R^3_1$ we use the algebra  $\D$ of the double numbers, which is determined in the following way:
$\D =\{a+\jj b  : \  a,b \in \R ,\ \jj^2=1 \}$,  where $\jj$ commutes with the elements of $\R$.
For the element  $z=a+\jj b$ of $\D$ we have $|z|^2=z\bar z = (a+\jj b)(a-\jj b) = a^2-b^2$.
This shows that $\D$ is the hyperbolic analogue of the algebra of complex numbers $\C$ and reflects the Lorentz geometry.
The algebra of double numbers is used essentially in paper \cite{Kanchev2020} in the study of the  Lorentz surfaces in $\R^4_2$.

Let $f(z)$ and $g(z)$ be two holomorphic functions defined in a domain of $\D$. Consider the
Weierstrass curve, defined by
$$
	\Psi(z)=\int_{z_0}^z \, 
	\left( -\frac12\, f(z)(1+g^2(z)),\frac{\jj}{2}\, f(z)(1-g^2(z)),f(z)g(z) \right) \, dz\ .
$$
 The real and the ''imaginary'' parts ${\bf x}(u,v)$ and ${\bf y}(u,v)$ define 
two  minimal  timelike surfaces of Gauss curvature $K<0$ and $-K$, respectively. 

\vskip 2mm
For example,  with the functions $f(z)=1,\ g(z)=z$ we obtain 
the classical Enneper surface of negative Gauss curvature 
$$
	{\bf x}(u,v)=\left( -\frac u6(u^2+3v^2+3), -\frac v6( 3u^2+v^2-3),\frac12(u^2+v^2) \right) 
$$
as the real part of the Weierstrass curve and the classical Enneper
surface of positive Gauss curvature
$$
	{\bf y}(u,v)=\left( -\frac v6(3u^2+v^2+3), -\frac u6 (u^2+3v^2-3),uv \right) 
$$
as the "imaginary"  part of the Weierstrass curve.

\begin{figure}[hbtp]
\begin{minipage}[b]{3.5in}
\centering
\includegraphics[width=0.7\textwidth]{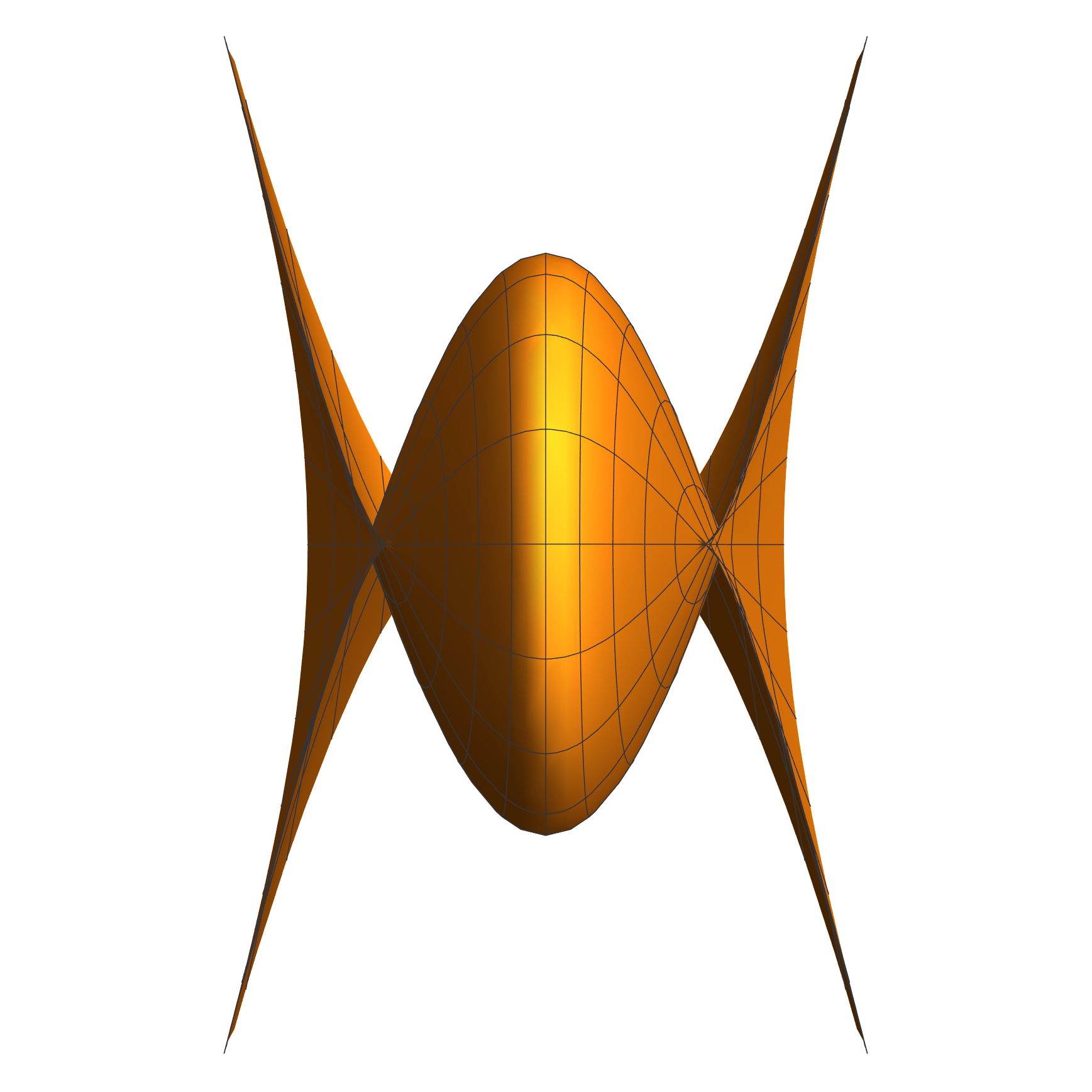}
 \centerline{\footnotesize Enneper surface with $K<0$}
  \end{minipage}
\begin{minipage}[b]{2.5in}
\centering
\includegraphics[width=1.1\textwidth]{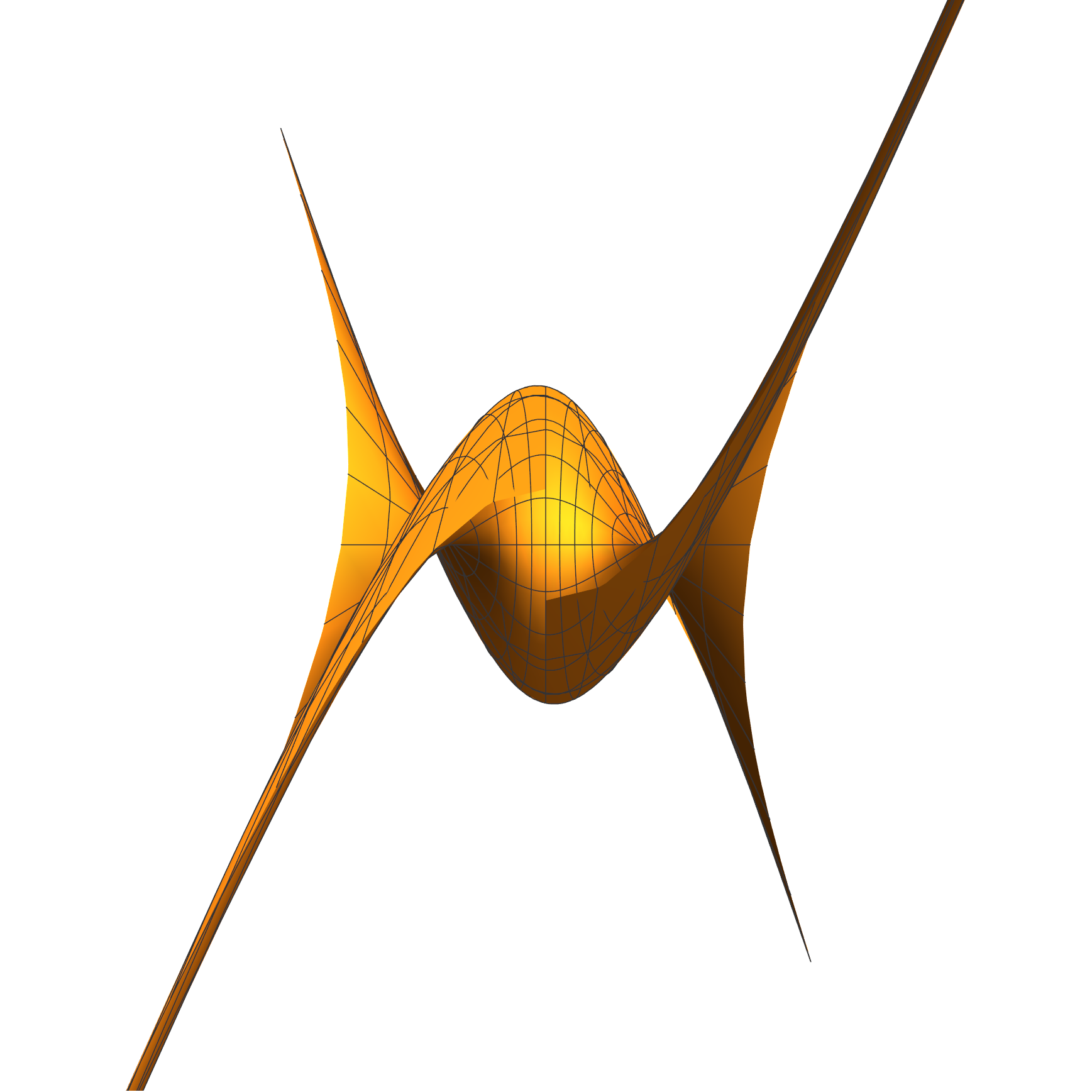}
    \centerline{\footnotesize Enneper surface with $K>0$}
    \end{minipage}
\caption{Enneper surfaces}\label{generating}
\end{figure}

Conversely, every minimal  timelike surface 
can be obtained  at least locally in this way. Note, however, that a minimal  timelike 
surface can be generated via the Weierstrass formula by different pairs of holomorphic 
functions on the algebra of double numbers. 

\vskip 2mm
In his study of  minimal timelike surfaces in $\mathbb R^3_1$, G. Ganchev \cite{G-time-like} 
specialized the Weierstrass formula by introducing special isothermal parameters, called 
 {\it canonical}. If the surface  
is parametrized with respect to canonical parameters,  then the coefficients of the first and second fundamental forms are
\begin{equation} \label{eq:principal}
	\begin{array}{l}
		\ds -E=G=\frac1{\sqrt{-K}} >0,  \qquad F=0, \vspace{1mm}\\
		\ds L=-1,  \qquad M=0, \qquad  N=-1,
	\end{array}
\end{equation}
in the case of minimal surfaces of negative Gauss curvature $K$, and
\begin{equation} \label{eq:asymptotic}
	\begin{array}{l}
		\ds E=-G=\frac1{\sqrt{K}} >0, \qquad F=0,\\
		\ds L=0,  \qquad M=1, \qquad  N=0,
	\end{array}
\end{equation}
in the case of minimal surfaces of positive Gauss curvature $K$.
Note that, because of (\ref{eq:principal}) the parametric lines are principal in the case 
of minimal surfaces with $K<0$. Analogously, because of (\ref{eq:asymptotic}) the parametric lines 
are asymptotic in the case of minimal surfaces with $K>0$.
The idea of Ganchev leads to the special Weierstrass curve 
\begin{equation} \label{eq:2.2}
	\Phi(z)=\int_{z_0}^z \, 
	\left( - \frac12\,\frac{1+g^2(z)}{g'(z)}, \frac{\jj}{2}\, \frac{1-g^2(z)}{g'(z)},\frac{g(z)}{g'(z)} \right) \, dz \ .  
\end{equation}
The real (resp. the ''imaginary'') part of this  curve is a minimal timelike 
surface with canonical  parametrization and negative (resp. positive) 
Gauss curvature. We shall use also the following theorem:     

\vspace{0.2 cm}
{\bf Theorem A \cite{G-time-like}.}  {\it If a timelike surface in $\mathbb R^3_1$ with non-vanishing Gauss curvature
is parametrized by canonical  parameters $(u,v)$,
then the Gauss curvature $K$ satisfies the equation
$$
	(\ln\sqrt{-K})_{uu}-(\ln\sqrt{-K})_{vv} =2\sqrt{-K}, 
$$  
in the case $K<0$, and
$$
	(\ln\sqrt{K})_{uu}-(\ln\sqrt{K})_{vv} =2\sqrt{K},
$$ 
in the case  $K>0$.
Conversely, for any solution $K(u,v)$ of any of these equations there exists
a {\bf unique} (up to position in the space) minimal timelike surface 
of  Gauss curvature $K(u,v)$, $(u,v)$ being canonical  parameters. }

\vspace{0.2 cm}
The canonical  parameters $(u,v)$ are determined 
uniquely up to the following transformations \cite{G-time-like}:
\begin{equation} \label{eq:2.3}
	\begin{array}{l}
	\vspace{2mm}
		u=\varepsilon\bar u+A\\
		\vspace{2mm}
		v=\varepsilon\bar v+B 
	\end{array} \qquad \varepsilon=\pm 1\ , \ A=const, \ B=const. 
\end{equation}

The idea of canonical parameters is further developed for the class of timelike surfaces
in  $\mathbb R^n_1$, see \cite{G-K}.

\setcounter{equation}{0}
\section{Transformation of the isothermal parameters to canonical ones}

Suppose the surface $S$ is defined as the real part of the Weierstrass  curve
\begin{equation} \label{eq:3.1}
   	\Psi(z)=\int_{z_0}^z \, 
	\left( -\frac12\, f(z)(1+g^2(z)),\frac{\jj}{2}\, f(z)(1-g^2(z)),f(z)g(z) \right) \, dz \ . 
\end{equation}
We look for a transformation $z=z(w)$ such that the curve $\Phi(w) = \Psi(z(w))$ has the form
$$
	\Phi(w)=\int_{w_0}^w \, 
	\left( - \frac12\,\frac{1+\tilde g^2(w)}{\tilde g'(w)}, 
	 \frac{\jj}{2}\,\frac{(1-\tilde g^2(w))}{\tilde g'(w)},\frac{\tilde g(w)}{\tilde g'(w)} \right) \, dw
$$
for some holomorphic function $\tilde g(w)$. The real part of this  curve will be 
a canonical  representation of the given surface $S$. The equality \ $\Psi(z(w))=\Phi(w)$ \
implies $\Psi'(z(w))z'(w)=\Phi'(w)$. Hence, it is easy to derive
\begin{equation} \label{eq:3.2}
	 f(z(w))z'(w)=\frac{1}{\tilde g'(w)} \ ,   \qquad\qquad  g(z(w))=\tilde g(w)  \ .  
\end{equation}
The last equality implies
$$
	\tilde g'(w)=g'(z(w))z'(w)
$$
and using the first equality of (\ref{eq:3.2}) we obtain
\begin{equation} \label{eq:3.3}
	(z'(w))^2=\frac1{f(z(w))g'(z(w))} \ .  
\end{equation}
Now, we know also the function $\tilde g(w)=g(z(w))$  that generates the surface in canonical  parameters.

Similar considerations can be done in the case the surface is defined as the  ''imaginary'' part of the Weierstrass curve determined by \eqref{eq:3.1}. 

So, we can state the following result:

\begin{thm}\label{T:3.1} 
Let the minimal timelike surface $S$ be defined by the real or ''imaginary''
part of (\ref{eq:3.1}). Any solution to differential equation (\ref{eq:3.3})
defines a transformation of the isothermal parameters of $S$ to canonical ones. 
Moreover, the function $\tilde g(w)$ that defines $S$ via formula (\ref{eq:2.2}) is given by $\tilde g(w)=g(z(w))$.
\end{thm}

As a consequence, we may obtain also   relations  (\ref{eq:2.3}) between two different pairs of canonical parameters.

\vspace{0.2cm}
As an application of Theorem 3.1, consider the minimal surfaces $S$ (of negative Gauss curvature) generated by
the functions
\begin{equation*} \label{eq:3.4}
	f(z) = a \ , \qquad g(z)=bz+c
\end{equation*}
via the Weierstrass formula, $a,b,c$ being double numbers, $a \neq 0$, $b \neq 0$, $ab\neq 0$. Equation (\ref{eq:3.3}) takes the form
$$
	(z'(w))^2=\frac{1}{ab}
$$
and has the following solution
$$
	z(w)=\pm\frac{w}{\sqrt a\sqrt b} + const.
$$
According to (\ref{eq:2.3}) and Theorem 3.1,  we may replace $z$ in $g(z)$ with  \ $\displaystyle{\frac{z}{\sqrt a\sqrt b}-\frac cb}$ \ and we will obtain a parametrization of
the surface $S$ in canonical  parameters via  formula (\ref{eq:2.2}) and the function
$$
	\tilde g(z)=g\left(\frac{z}{\sqrt a\sqrt b}-\frac cb\right)=\frac{\sqrt b}{\sqrt a}\, z \ .
$$
The Gauss curvature is given by the following formula:
$$
	K=-\frac{16\left| \frac ba \right|^2}{\left(1-\left| \frac ba \right|(u^2-v^2)\right)^4} \ .
$$ 
This result shows that the surface $S_0$ of negative Gauss curvature generated via the 
Weierstrass formula by the functions
$$
	f(z) =\frac{|a|}{|b|} \ , \qquad g(z)=z,
$$
has the same Gauss curvature in canonical  parameters. Hence, due to Theorem A
we may identify $S$ with $S_0$. On the other hand, the Weierstrass formula implies 
that $S_0$ is homothetic to the standard timelike Enneper surface ($a=b=1$) with $K<0$. 

\vskip 2mm
So, as in the case of minimal surfaces in the Euclidean space \cite{C-M}, we have

\begin{cor}\label{T:3.3}  The minimal timelike surface  generated by the pair of linear functions 
	$f(z) = a, \, g(z)=bz+c$ via the Weierstrass formula coincides with the Enneper surface up to position in the space and homothety.
\end{cor}

\vspace{3mm}

\setcounter{equation}{0}
\section{Holomorphic functions generating a minimal timelike surface}

As we said before, a minimal timelike surface is generated by different pairs of 
holomorphic functions via the Weierstrass formula.  For example, the Enneper surface
of negative (resp. positive) Gauss curvature is the real (resp.  ''imaginary'') part
of the  curve defined via the Weierstrass formula by the pair of functions
\begin{equation} \label{eq:4.0a}
	f(z)=1;   \qquad  g(z)=z,
\end{equation}
but also by
\begin{equation} \label{eq:4.0b}
	f(z)=e^{z};   \qquad  g(z)=e^z,
\end{equation}
and, of course, by many others. So, the following natural question arises: under what conditions do two 
pairs of holomophic functions give rise to one and the same  minimal timelike
surface via the Weierstrass representation? It is not difficult
to prove the following:

\begin{prop}\label{T:4.1}
Suppose the pairs 
$	(\tilde f(z),\tilde g(z))$ and $ (f(w),g(w))$
generate two minimal timelike surfaces via the Weierstrass formula. Then, these 
surfaces coincide (up to translation) if and only if there exists a function $w=w(z)$, such that
$$
	\tilde g(z)=g(w(z)) 
\qquad {\rm and}  \qquad
	\tilde f(z)=f(w(z))w'(z)   \ .
$$
\end{prop}

For the  two pairs \eqref{eq:4.0a} and \eqref{eq:4.0b}, that generate the Enneper surface, the function  $w(z)=e^z$
transfers the first pair into the second one.

Similarly, the following question related to   formula (\ref{eq:2.2})  arises:  
what is the relation between  the functions that generate a  minimal timelike surface in canonical 
 parameters? A result in this direction is given by the following theorem.

\vspace{0.2cm}
\begin{thm}\label{T:4.3}  Let the holomorphic function $ g(z)$ (defined on a domain of $\D$) generate a minimal timelike surface $S$ in canonical  parameters, i.e. via  formula (\ref{eq:2.2}). 
Then, for an arbitrary real number $\varphi$ and an arbitrary double number $\alpha$, by the 
 transformations
\begin{equation} \label{eq:4.1}
  \tilde g(z)=\pm e^{\varphi j}\frac{\alpha+g(z)}{ 1+\bar \alpha g(z)};  
	\qquad\qquad
	\tilde g(z)=\pm e^{\varphi j}\frac1{f(z)},  
\end{equation}  
we obtain the same (up to position in the space) surface  in canonical  parameters. 
Conversely,  any function that generates  (up to position) the 
surface $S$ in canonical  parameters may be obtained in this way.
\end{thm}

{\bf Proof.} Let us consider the first transformation.
Denote by $S$ the minimal timelike surface of negative Gauss curvature, 
generated  via formula (\ref{eq:2.2}) by the function $g(z)$ and let  
$\Psi(z)$ be  the  corresponding  curve. Analogously, we define 
$\widetilde S$ and $\widetilde\Psi(z)$.

We may prove that $ S$ and $\widetilde S$ coincide (up to position) by a 
direct computation of their Gauss curvatures using the formula
$$
	K=-\frac{16|g'|^4}{(1-|g|^2)^4}       
$$
and applying Theorem A. Now we give another  proof, thus
clarifying the relation between $ S$, $\widetilde S$ and  transformation 
(\ref{eq:4.1}). 
We have
$$
	\Psi'(z)=\left(-  \frac{1+ g^2(z)}{2 g'(z)},
	         \jj \frac{1- g^2(z)}{2 g'(z)},\frac{ g(z)}{ g'(z)} \right)  \ ,
$$
$$
	\widetilde\Psi'(z)=\left(-  \frac{1+\widetilde g^2(z)}{2\widetilde g'(z)},
	         \jj \frac{1-\widetilde g^2(z)}{2\widetilde g'(z)},\frac{\widetilde g(z)}{\widetilde g'(z)} \right) \ .
$$

Let $\alpha=a+ \jj b$, $a,b\in \R$. Define the SO(1,2)-matrices
\begin{equation*}
	A=\left(\begin{array}{ccc}
	     \cosh\varphi  &   \sinh\varphi  &  0  \\
	     \sinh\varphi  &   \cosh\varphi  &  0  \\
	          0       &       0        &  1
	  \end{array}\right);   \qquad  \qquad
	  	B=\left(\begin{array}{ccc}
	     \frac{1+a^2+b^2}{1-a^2+b^2}  &  \frac{-2ab}{1-a^2+b^2}  &  \frac{-2a}{1-a^2+b^2}  \vspace{.1cm} \\ 
	     \frac{2ab}{1-a^2+b^2}  &  \frac{1-a^2-b^2}{1-a^2+b^2}   & \frac{-2b}{1-a^2+b^2}   \vspace{.1cm} \\
	     \frac{-2a}{1-a^2+b^2}   &  \frac{2b}{1-a^2+b^2}          &  \frac{1+a^2-b^2}{1-a^2+b^2}
	  \end{array}\right).
\end{equation*}
A straightforward verification shows that 
$$
	A\,B\,\Psi'(z)=\widetilde\Psi'(z)  \ .
$$
The last equality implies that up to translation
$$
	A\,B\,{\bf x}(u,v)=\tilde{\bf x}(u,v)  \ .
$$
Hence, the considered transformation of type (\ref{eq:4.1}) of the function $g(z)$ corresponds to a motion
of the surface $S$. 

Conversely,  it is clear  that any surface that coincides (up to position) with $S$ may be 
obtained from $S$ using as above two $SO(1,2)$ matrices and a translation.

\hfill{\qed}

\vspace{0.2cm}
As an application of Theorem \ref{T:4.3}, we may prove that any polynomial 
minimal surfaces $S$ 	which has polynomial parametrization of degree 3 in
isothermal coordinates is (up to position and homothety) an Enneper surface. 
Namely, we have the following result:

\begin{thm}\label{T:4.4}  Let the minimal timelike surface $S$ of negative Gauss curvature has  
polynomial parametrization of degree 3 in isothermal parameters. Then, up to position in space 
and homothety,  $S$ is (a part of) the Enneper surface of negative curvature. 
\end{thm}  

{\bf Proof.} Suppose that the surface is defined by
$$
	S\ : \ \ {\bf x}={\bf x}(u,v).
$$
Similarly to the case of surfaces in the Euclidean space (see e.g. section 22.4 in \cite{G-A-S}),
$S$ is the real part of the Weierstrass curve obtained by substituting 
the double number variables $\displaystyle{\frac{z}{2}}$ and $\displaystyle{\frac{z}{2\jj}}$ formally in the places of the real variables $u, v$:
$$
	\Psi(z)=2{\bf x}\left(\frac z2,\frac z{2 \jj}\right)-{\bf x}(0,0).
$$
Using a translation (if necessary) we may assume that ${\bf x}(0,0)=0$. Since the curve $\Psi(z)$ is  a cubic polynomial, then
$$
	\Psi'(z)=(\phi_1(z),\phi_2(z),\phi_3(z))= 
	\left( -\frac{f(z)}2(1+g^2(z)),\frac{\jj}{2} f(z)(1-g^2(z)), f(z)g(z) \right) \, dz\ 
$$
for some functions $f(z)$ and $g(z)$ and the functions $\phi_i(z)$, $i =1, 2, 3$ are polynomials
of degrees at most 2. Moreover, at least one of them is of degree exactly 2. Hence, the same is
true for the following three functions
\begin{equation} \label{eq:3.5}
	f(z)=-\phi_1+ \jj \phi_2;  \qquad  f(z)g^2(z)=-\phi_1- \jj \phi_2;   \qquad   f(z)g(z)=\phi_3 \ .
\end{equation}
So, $f(z)$ is a polynomial of degree at most 2. 
From the third equality of \eqref{eq:3.5} we have $g=\displaystyle{\frac{\phi_3}{f}}$. Since $\phi_3$ and $f$
are polynomials, we may write $g(z)$ in the form 
$$
	g(z)=\frac{P(z)}{Q(z)} \ ,
$$  
where the polynomials $P(z)$ and $Q(z)$ have no common zeros. If  we assume that $Q(z)$ is a constant, then
$g(z)$ is a polynomial and having in mind \eqref{eq:3.5}, we get $f(z)=const$, $g(z)=cz+d$.
If we assume that $Q(z)$ is not a constant, then from the second equality of  \eqref{eq:3.5}, we get 
$$
	f(z)\frac{P^2(z)}{Q^2(z)}=-\phi_1-j\phi_2,
$$
which is a polynomial and hence
$$
	f(z)=\pm (az+b)^2 =\pm Q^2(z)\ .
$$
Up to symmetry of the surface, we assume that
$$
	f(z)=(az+b)^2, \qquad Q(z)=az+b \ .
$$
Now, $\psi_3=f(z)g(z)=(az+b)P(z)$ and since it is of degree at most 2, then $P(z)$ is of degree at most 1,
i.e. $P(z)=cz+d$.

Hence, we conclude that, up to homothety of the surface, $f(z)$ and $g(z)$ have the form
$$
	f(z)=(az+b)^2;   \qquad g(z)=\frac{cz+d}{az+b}  \ .
$$ 

On the other hand, the Enneper surface is generated in canonical parameters by the pair 
of functions $f_1(z)=1$, $g_1(z)=z$ and due to Theorem \ref{T:4.3} also by the functions
$$
	 g_2(z)=e^{\varphi \jj}\frac{\alpha+z}{1+\bar\alpha z};
		\qquad\qquad
	 f_2(z)=\frac1{g_2'(z)}=e^{-\varphi \jj}\frac{(1+\bar\alpha z)^2}{|\alpha|^2-1}.
$$
We can change  the parameter $z$ by	
$$
	\frac{(a\alpha e^{\varphi \jj}-c)z+\alpha b e^{\varphi \jj}-d}{(\bar\alpha c-ae^{\varphi \jj})z+\bar\alpha d-be^{\varphi \jj}}
$$
and then the generating functions take the form:
$$
	f_3(z)=\frac{e^{\varphi \jj}(1-|\alpha|^2)(az+b)^2}{\Big( (\bar\alpha c-ae^{\varphi \jj})z+\bar\alpha d-be^{\varphi \jj} \Big)^2};
	\qquad\qquad
	g_3(z)=\frac{cz+d}{az+b}.
$$
Of course, in the last expressions the parameters are not canonical. Note that we may choose arbitrary $\varphi$ and $\alpha$, so we put
$$
	e^{\varphi \jj}=\frac{\frac{c^2}{(bc-ad)^2}}{\left|\frac{c^2}{(bc-ad)^2}\right|}; 
	\qquad\qquad 
	\bar\alpha=\frac{ae^{\varphi \jj}}{c} \ .
$$
Then, $f_3(z)$ becomes proportional (with real coefficient) to $f(z)$, thus proving the assertion.
Note that $bc-ad$ can not be zero (if we assume that $bc - ad =  0$, then the surface is planar which is 
not our case).  

\hfill{\qed}

\begin{remark} The same proposition holds for minimal timelike surfaces of positive Gauss curvature.
\end{remark}

\vskip 5mm \textbf{Acknowledgments:}
The  authors are partially supported by the National Science Fund,
Ministry of Education and Science of Bulgaria under contract KP-06-N52/3.


\begin{thebibliography}{9}


\bibitem{Ch-D-Mag}
R.~M.~B. Chaves,  M.~P. Dussan, M. Magid,
Bj\"orling problem for timelike surfaces in the Lorentz-Minkowski space. Math. Anal. Appl. \textbf{377} (2), (2011), 481--494.
 doi:\,10.1016/j.jmaa.2010.10.076.

\bibitem{Cintra-Onnis}
 A.~A. Cintra, I.~I. Onnis, Enneper representation of minimal surfaces in the three-dimensional Lorentz-Minkowski space.
Annali di Matematica, \textbf{197} (1), (2018), 21--39.
doi:\,10.1007/s10231-017-0666\nobreakdash-z.

\bibitem{C-M}
C. Cos\'in, J. Monterde, B\'ezier Surfaces of Minimal Area. Proc.
Int. Workshop of Computer Graphics and Geom. Modelig. Lecture Notes of
Computer Science,  72--81, Springer-Verlag, 2002. 

\bibitem{Erd}
 S. Erdem, Harmonic maps of Lorentz surfaces, quadratic differentials and paraholomorphicity. Beitr\"age Algebra Geom. 38 (1) (1997) 19--32.

\bibitem{Fu-In}
 A. Fujioka, J. Inoguchi, Timelike surfaces with harmonic inverse mean curvature, in: Surveys on Geometry and Integrable Systems, in: Adv. Stud. Pure
Math., vol. 51, Math. Soc. of Japan, Tokyo, 2008, pp. 113--141.


\bibitem{G-2008} G. Ganchev, Canonical Weierstrass representation of minimal surfaces in Euclidean space.
ArXiv: 0802.2374.

\bibitem{G-time-like}
G. Ganchev, Canonical Representations of Minimal Time-like Surfaces in Minkowski Space and Explicit Solving of Their Natural PDE. (preprint) 

\bibitem{G-K} 
G. Ganchev, K. Kanchev, Canonical coordinates on minimal time-like surfaces in the $n$-dimensional 
Minkowski space. Serdica Math. J. \textbf{45} (2019), 341--372.

\bibitem{Kanchev2020}
G. Ganchev, K. Kanchev,
Canonical coordinates and natural equations for minimal time-like
  surfaces in $\mathbb{R}^4_2$.  Kodai Mathematical Journal, 43 (3)  2020, 524--572.


\bibitem{G-M} 
G. Ganchev, V. Mihova, On the invariant theory of Weingarten surfaces in Euclidean space.  
J. Phys. A, {\bf 43}, 40(2010),405210, 27 pp.

\bibitem{G-A-S} 
A. Gray, E. Abbena, S. Salamon, Modern Differential Geometry of Curves and Surfaces with Mathematica, Third Edition. Chapman and Hall/CRC, 2006.

\bibitem{In-Tod}
 J. Inoguchi, M. Toda, Timelike minimal surfaces via loop groups, Acta Appl. Math. 83 (3) (2004) 313--355.

\bibitem{OK} 
O. Kassabov, Transition to canonical principal parameters on minimal surfaces. Comput.
Aided Geom. Design. {\bf 31} (2014), 441--450.

\bibitem{Kond}
J. Konderak,  A Weierstrass representation theorem for Lorentz surfaces. Complex Var. Theory Appl., \textbf{50}  (2005), 319--332.
doi:\,10.1080/02781070500032895.

\bibitem{Lira-Melo-Merc}
J.~H. Lira, M. Melo,  F. Mercuri,  A Weierstrass representation for minimal surfaces in 3-dimensional manifolds.
Results. Math., \textbf{60}, (2011), 311--323. doi:\,10.1007/s00025-011-0169-y.


\bibitem{Mag}
M. Magid, Timelike surfaces in Lorentz 3-space with prescribed mean curvature and Gauss map, Hokkaido Math. J. 19 (1991), 447--464. DOI: 10.14492/hokmj/1381413979

\bibitem{Ship-Pack}
B.~A. Shipman, P.~D. Shipman,  D. Packard, Generalized Weierstrass-Enneper representations of Euclidean, spacelike, and timelike surfaces: a unified Lie-algebraic formulation.
 J. Geom., \textbf{108} (2), (2017), 545--563. doi:\,10.1007/s00022-016-0358-7.

\bibitem{Weier} 
K.T.W. Weierstrass, Untersuchung \"uber die Fl\"achen, deren mittlere Kr\"ummung \"uberall
gleich null ist, Monatsberichte der Berliner Akademie, 1866, 612--625.

\end{thebibliography}
\end{document}